\documentclass[11pt,leqno,fleqn]{article}
\usepackage{amsmath,amssymb,hyperref}
\usepackage{amsfonts}
\newcommand{\NIET}[1]{}
\newcommand{\CC}{{\cal C}}
\newcommand{\oR}{{\mathbb{R}}}
\newcounter{sectie}
\newcounter{subsectie}
\newcommand{\sectz}[1]{\refstepcounter{sectie}\setcounter{subsectie}{0}
\section*{\boldmath \thesectie. #1}%
}
\newcommand{\oZ}{{\mathbb{Z}}}
\newcommand{\dez}[1]{\dyz{\raggedright$\displaystyle #1 $}}
\newcommand{\dyz}[1]{%
\refstepcounter{equation}%
\begin{list}{}{
\topsep 3mm
\leftmargin 18mm
\rightmargin 0cm
\itemsep 0mm
\listparindent 0mm
\parsep 0mm
\itemsep 0mm
\labelsep 0mm
\labelwidth 18mm
}%
\item[\rm (\theequation)\hfill]
#1
\end{list}%
}
\newcommand{\dy}[2]{%
\refstepcounter{equation}%
\label{#1}%
\begin{list}{}{
\topsep 3mm
\leftmargin 18mm
\rightmargin 0cm
\itemsep 0mm
\listparindent 0mm
\parsep 0mm
\itemsep 0mm
\labelsep 0mm
\labelwidth 18mm
}%
\item[\rm (\theequation)\hfill]
#2
\end{list}%
}
\newcommand{\dps}{\displaystyle}
\newcounter{bewering}
\newcommand{\propz}[1]{\refstepcounter{bewering}\vspace{4mm}\noindent{\bf Proposition \thebewering.}{\it #1}}
\newcommand{\pf}{\vspace{3mm}\noindent{\bf Proof.}\ }
\newcommand{\bx}{\hspace*{\fill} \hbox{\hskip 1pt \vrule width 4pt height 8pt depth 1.5pt \hskip 1pt}

\addvspace{4mm}}
\newcommand{\rf}[1]{{\rm (\ref{#1})}}
\newcommand{\T}{^{\sf T}}
\newcommand{\rondje}[1]{\accent23 #1}
\newcommand{\subsectz}[1]{\refstepcounter{subsectie}
\subsection*{\boldmath \thesectie.\thesubsectie. #1}%
}
\newcommand{\sect}[2]{\refstepcounter{sectie}\setcounter{subsectie}{0}
\section*{\boldmath \thesectie. #2}%
\label{#1}}
\newcommand{\dyy}[2]{\dy{#1}{\raggedright$\dps#2$}}
\newcommand{\de}[2]{\dy{#1}{\raggedright$\displaystyle #2 $}}
\newcommand{\subsect}[2]{\refstepcounter{subsectie}
\subsection*{\boldmath \thesectie.\thesubsectie. #2}%
\label{#1}}
\newcommand{\dyyz}[1]{\dyz{\raggedright$\dps#1$}}
\newcommand{\sgn}{\text{\rm sgn}}
\newcommand{\bm}[1]{\text{\boldmath$#1$}}
\newcommand{\OO}{{\cal O}}
\newcommand{\prop}[2]{\refstepcounter{bewering}\vspace{4mm}\noindent{\bf Proposition \thebewering.}\label{#1}{\it #2}}
\newcommand{\Sym}{{\text{\rm Sym}}}
\newcommand{\kies}[2]{\mbox{${{#1}\choose{#2}}$}}
\renewcommand{\part}{\text{\rm part}}
\newcounter{hulpstelling}
\newcommand{\lemmaz}[1]{\refstepcounter{hulpstelling}\vspace{4mm}\noindent{\bf Lemma \thehulpstelling.}{\it #1}}
\newcommand{\height}{\text{\rm height}}
\begin{document}

\begin{center}
{\large\bf SEMIDEFINITE BOUNDS FOR NONBINARY CODES BASED ON QUADRUPLES

}
\vspace{4mm}
Bart Litjens\footnote{ Korteweg-De Vries Institute for Mathematics,
University of Amsterdam.
The research leading to these results has received funding from the European Research Council
under the European Union's Seventh Framework Programme (FP7/2007-2013) / ERC grant agreement
n$\mbox{}^{\circ}$ 339109.},
Sven Polak$\mbox{}^1$,
Alexander Schrijver$\mbox{}^1$

\end{center}

\noindent
{\small
{\bf Abstract.}
For nonnegative integers $q,n,d$, let $A_q(n,d)$ denote the maximum cardinality of a code of length $n$
over an alphabet $[q]$ with $q$ letters and with minimum distance at least $d$.
We consider the following upper bound on $A_q(n,d)$.
For any $k$, let $\CC_k$ be the collection of codes of cardinality at most $k$.
Then $A_q(n,d)$ is at most the maximum value of $\sum_{v\in[q]^n}x(\{v\})$, where
$x$ is a function $\CC_4\to\oR_+$ such that $x(\emptyset)=1$ and $x(C)=0$ if $C$ has minimum distance less than $d$,
and such that the $\CC_2\times\CC_2$ matrix $(x(C\cup C'))_{C,C'\in\CC_2}$ is positive semidefinite.
By the symmetry of the problem, we can apply representation theory to reduce the problem to a
semidefinite programming problem with order bounded by a polynomial in $n$.
It yields the new upper bounds
$A_4(6,3)\leq 176$,
$A_4(7,4)\leq 155$,
$A_5(7,4)\leq 489$, and
$A_5(7,5)\leq 87$.

}

\medskip
\noindent
{\bf Key words:} code, nonbinary code, upper bounds, semidefinite programming, Delsarte\\
{\bf MSC 2010:} 94B65, 05E10, 90C22, 20C30

\sectz{Introduction}

Let $\oZ_+$ denote the set of nonnegative integers, and denote $[m]=\{1,\ldots,m\}$, for any $m\in\oZ_+$.
Fixing $n,q\in\oZ_+$, a {\em code} is a subset of $[q]^n$.
So $[q]$ serves as the alphabet and $n$ as the word length.
We will assume throughout that $q\geq 2$.
(If you prefer $\{0,1,\ldots,q-1\}$ as alphabet, take the letters mod $q$.)
While this paper is mainly meant to handle the case $q\geq 3$, the results also hold for $q=2$.

For $v,w\in[q]^n$, the {\em (Hamming) distance} $d_H(v,w)$ is equal to the number of $i\in[n]$ with
$v_i\neq w_i$.
The {\em minimum distance} of a code $C$ is the minimum of $d_H(v,w)$ taken over distinct $v,w\in C$.
Then $A_q(n,d)$ denotes the maximum cardinality of a code with minimum distance at least $d$.
We will study the following upper bound on $A_q(n,d)$, sharpening
Delsarte's classical linear programming bound [4].

For $k\in\oZ_+$,
let $\CC_k$ be the collection of subsets $C$ of $[q]^n$ with $|C|\leq k$.
For each $x:\CC_4\to\oR$ define the $\CC_2\times\CC_2$ matrix $M(x)$ by
\dez{
M(x)_{C,C'}:=x(C\cup C')
}
for $C,C'\in \CC_2$.
Then define
\dy{8no15a}{
$B_q(n,d):=\dps\max_x\sum_{w\in[q]^n}x(\{w\})$,
where $x:\CC_4\to\oR_+$ satisfies\\
(i) $x(\emptyset)=1$,\\
(ii) $x(C)=0$ if the minimum distance of $C$ is less than $d$,\\
(iii) $M(x)$ is positive semidefinite.
}

\propz{
$A_q(n,d)\leq B_q(n,d)$.
}

\pf
Let $D\subseteq[q]^n$ have minimum distance at least $d$ and satisfy $|D|=A_q(n,d)$.
Define $x:\CC_4\to\oR$ by $x(C)=1$ if $C\subseteq D$ and $x(C)=0$ otherwise.
Then $x$ satisfies the conditions: (iii) follows from the fact that for this $x$ one has $M(x)_{C,C'}=x(C)x(C')$
for all $C,C'\in\CC_2$.
Moreover, $\sum_{w\in[q]^n}x(\{w\})=|D|=A_q(n,d)$.
\bx

The optimization problem \rf{8no15a} is huge, but, with methods from representation theory,
can be reduced to a size bounded by a polynomial in $n$, with entries (i.e.,
coefficients) being polynomials in $q$.
This makes it possible to solve \rf{8no15a} by semidefinite programming for some moderate values
of $n$, $d$, and $q$, leading to improvements of best known upper bounds for $A_q(n,d)$.

To explain the reduction, let $H$ be the wreath product $S_q^n\rtimes S_n$.
For each $k$, the group $H$ acts naturally on $\CC_k$, maintaining minimum distances and cardinalities
of elements of $\CC_k$ (being codes).
Then we can assume that $x$ is invariant under the $H$-action on $\CC_4$.
That is, we can assume that $x(C)=x(D)$ whenever $C,D\in\CC_2$ and $D=g\cdot C$ for some $g\in H$.
Indeed, \rf{8no15a}(i)(ii)(iii) are maintained under replacing $x$ by $g\cdot x$.
(Note that $M(g\cdot x)$ is obtained from $M(x)$ by simultaneously permuting rows and columns.)
Moreover, the objective function does not change by this action.
Hence the optimum $x$ can be replaced by the average of all $g\cdot x$ (over all $g\in  H$),
by the convexity of the set of positive semidefinite matrices.
This makes the optimum solution $H$-invariant.

Let $\Omega$ be the set of $H$-orbits on $\CC_4$.
Note that $\Omega$ is bounded by a polynomial in $n$ (independently of $q$).
As there exists an $H$-invariant optimum solution, we can replace, for each $\omega\in\Omega$ and $C\in\omega$, each variable $x(C)$ by a variable $y(\omega)$.
In this way we obtain $M(y)$.

Then $M(y)$ is invariant under the action of $H$ on its rows and columns,
induced from the action of $H$ on $\CC_2$.
Hence $M(y)$ can be block-diagonalized by $M(y)\mapsto U\T M(y) U$, where $U$ is a
matrix independent of $y$.
The entries in each block are linear functions of the variables $y(\omega)$.
There are several equal (or equivalent) blocks.
Taking one block from each such class gives a matrix of order polynomial in $n$ with numbers
that are polynomials in $q$.
The issue crucial for us is that the original matrix $M(y)$ is positive semidefinite if and only if
each of the blocks is positive semidefinite.

In this paper we will describe the blocks that reduce the problem.
With this,
we found the following improvements on the known bounds for $A_q(n,d)$, with thanks to Hans D. Mittelmann
for his help in solving the larger-sized problems

\begin{center}
\begin{tabular}{|c|c|c||c|c|c|}
\hline
&&&           &          &best\\
&&&best       &          &upper\\
&&&lower      &{\bf new}       &bound\\
&&&bound      &{\bf upper}      &previously\\
$q$&$n$&$d$&known      &{\bf bound}      &known\\
\hline
\hline
4&6&3&164&{\bf 176}&179\\
\hline
4&7&4&128&{\bf 155}&169\\
\hline
\hline
5&7&4&250&{\bf 489}&545\\
\hline
5&7&5&\hspace*{5.5pt}53&\hspace*{5.5pt}{\bf 87}&108\\
\hline
\end{tabular}
\end{center}
The best upper bound 179 previously known for $A_4(6,3)$ is the Delsarte bound [4];
the other three best upper bounds previously known were given by Gijswijt, Schrijver, and Tanaka [7].
We refer to the most invaluable tables maintained by Andries Brouwer [3] with the best
known lower and upper bounds for the size of error-correcting codes (see also
Bogdanova, Brouwer, Kapralov, and \"Osterg\rondje ard [1] and Bogdanova and \"Osterg\rondje ard [2] for studies of bounds for codes over alphabets of
size $q=4$ and $q=5$, respectively).

\subsectz{Comparison with earlier bounds}

The bound $B_q(n,d)$ described above is a sharpening of Delsarte's classical linear programming bound [4].
The value of the Delsarte bound is equal to our bound after replacing $\CC_4$ and $\CC_2$ by
$\CC_2$ and $\CC_1$, respectively, which generally yields
a less strict bound.

We can add to \rf{8no15a} the condition that, for each $D\in\CC_4$, the $S(D)\times S(D)$ matrix
\dy{20de15a}{
$(x(C\cup C'))_{C,C'\in S(D)}$ is positive semidefinite,
}
where $S(D):=\{C\in\CC_4\mid C\supseteq D, |D|+2|C\setminus D|\leq 4\}$.
(So (iii) in \rf{8no15a} is the case $D=\emptyset$.)
Also the addition of \rf{20de15a} allows a reduction of the optimization problem to polynomial
size as above.
For $q=2$ we obtain in this way the bound given by Gijswijt, Mittelmann, and Schrijver [6].
Our present description gives a more conceptual and representation-theoretic approach to the method of
[6].
Some preliminary computations suggest that adding condition \rf{20de15a} might be superfluous in the
optimization problem, but this needs further investigation.

A bound intermediate to the Delsarte bound and the currently investigated bound is based on
considering functions $x:\CC_3\to\oR_+$ and the related matrices ---
see 
Schrijver [9] for binary codes and
Gijswijt, Schrijver, and Tanaka [7] for nonbinary codes.

\sect{16de15j}{Preliminaries on representation theory}

We assume some familiarity with classical representation theory, in particular of the symmetric group $S_n$
and of finite groups in general.
In this section we give a brief review, also to settle some notation and terminology.
We refer to Sagan [8] for background.
All groups considered are finite, which allows us to keep decompositions and reductions real-valued.

A group $G$ {\em acts} on a set $X$ if there is a group homomorphism
$G\to S_X$, where $S_X$ is the set of bijections $X\to X$.
The image of $g\in G$ in $S_X$ is indicated by $g\cdot\mbox{}$.
If $X$ is a linear space, the bijections are assumed to be linear functions.
The action of $G$ on a set $X$ induces an action of $G$ on the linear
space $\oR^X$, by $(g\cdot f)(x):=f(g^{-1}\cdot x)$
for all $g\in G$, $f\in\oR^X$, and $x\in X$.
If any group $G$ acts on $X$, then $X^G$ denotes the set of elements of $X$
invariant under the action of $G$.

Let $m\in\oZ_+$ and let $G$ be a finite group acting on $V=\oR^m$.
Then $V$ can be decomposed uniquely as direct sum of the $G$-isotypical components
$V_1,\ldots,V_k$.
Next, each $V_i$ is a direct sum $V_{i,1}\oplus\cdots\oplus V_{i,m_i}$
of mutually $G$-isomorphic, irreducible $G$-modules.
(This decomposition is generally not unique.)
For each $i\leq k$ and $j\leq m_i$, choose a nonzero
$u_{i,j}\in V_{i,j}$ such that for each $i$ and all
$j,j'\leq m_i$ there exists a $G$-isomorphism $V_{i,j}\to V_{i,j'}$ bringing $u_{i,j}$ to $u_{i,j'}$.
For each $i\leq k$,
let $U_i$ be the matrix $[u_{i,1},\ldots,u_{i,m_i}]$, considering the $u_{i,j}$ as columns.
We call any matrix set $\{U_1,\ldots,U_k\}$ that can be obtained in this way {\em representative} for
the action of $G$ on $\oR^m$.
It has the property that the function
\dyy{11de15a}{
\Phi:(\oR^{m\times m})^G\to\bigoplus_{i=1}^k\oR^{m_i\times m_i}
\text{~~~~with~~~~}
\Phi(X):=\bigoplus_{i=1}^kU_i\T XU_i
}
for $X\in(\oR^{m\times m})^G$ is bijective.
So $\sum_im_i^2$ is equal to the dimension of $(\oR^{m\times m})^G$ (and hence can be considerably
smaller than $m$).
Another important property of a representative matrix set is that
$X$ is positive semidefinite if and only if $\Phi(X)$ is positive semidefinite.
Moreover, by construction, as $V_{i,j}=\oR G\cdot u_{i,j}$ for all $i,j$,
\de{28de15ster}{
V=\bigoplus_{i=1}^k\bigoplus_{j=1}^{m_i}\oR G\cdot u_{i,j}.
}

It will turn out to be convenient to consider the columns of the matrices $U_i$ as
elements of the dual space $(\oR^m)^*$ (by taking the standard inner product).
Then each $U_i$ is an ordered set of linear functions on $\oR^m$.
(The order plays a role in describing a representative matrix set for the action of the
wreath product $G\rtimes S_n$ on $V^{\otimes n}$.)

\subsect{16de15h}{A representative set for the action of $S_n$ on $V^{\otimes n}$}

Classical representation theory of the symmetric group yields
a representative set for the natural action of $S_n$
on $V^{\otimes n}$, where $V$ is a finite-dimensional vector space, which we
will describe now.

For $n\in\oZ_+$, $\lambda\vdash n$ means that $\lambda$ is equal to
$(\lambda_1,\ldots,\lambda_t)$ for some $t$,
with $\lambda_1\geq\cdots\geq\lambda_t>0$ integer and $\lambda_1+\cdots+\lambda_t=n$.
The number $t$ is called the {\em height} of $\lambda$.
The {\em Young shape} $Y(\lambda)$ of $\lambda$ is the set
\dyyz{
Y(\lambda):=\{(i,j)\in\oZ_+^2\mid 1\leq j\leq t, 1\leq i\leq\lambda_j\}.
}
For any $j_0\leq t$, the set of elements $(i,j_0)$ in $Y(\lambda)$ is called the
{\em $j_0$-th row} of $Y(\lambda)$.
Let $R_{\lambda}$ be the group of permutations $\pi$ of $Y(\lambda)$ with
$\pi(Z)=Z$ for each row $Z$ of $Y(\lambda)$.
For any $i_0\leq \lambda_1$, the set of elements $(i_0,j)$ in $Y(\lambda)$ is called the
{\em $i_0$-th column} of $Y(\lambda)$.
Let $C_{\lambda}$ be the group of permutations $\pi$ of $Y(\lambda)$ with
$\pi(Z)=Z$ for each column $Z$ of $Y(\lambda)$.

A {\em $\lambda$-tableau} is a function $\tau:Y(\lambda)\to\oZ_+$.
We put $\tau\sim\tau'$ for $\lambda$-tableaux $\tau,\tau'$ if $\tau'=\tau r$
for some $r\in R_{\lambda}$.
A $\lambda$-tableau is {\em semistandard} if
in each row the entries are nondecreasing and in each column the entries are
increasing.
Let $T_{\lambda,m}$ denote the collection of semistandard $\lambda$-tableaux
with entries in $[m]$.
Note that $T_{\lambda,m}\neq\emptyset$ if and only if $\lambda$ has height at most $m$.

Let $B=(B(1),\ldots,B(m))$ be an ordered basis of $V^*$.
For $\tau\in T_{\lambda,m}$, define the following element of $(V^*)^{\otimes n}$:
\de{14de15f}{
u_{\tau,B}:=\sum_{\tau'\sim\tau}\sum_{c\in C_{\lambda}}\sgn(c)\bigotimes_{y\in Y(\lambda)}B(\tau'c(y)),
}
where we order the Young shape $Y(\lambda)$ by concatenating its rows.
Then the matrix set
\dyy{16se15a}{
\{~~
[
u_{\tau,B}\mid \tau\in T_{\lambda,m}]
~~\mid~~
\lambda\vdash n
\}
}
is representative for the natural action of $S_n$ on $V^{\otimes n}$.

\sect{16de15i}{Reduction of the optimization problem}

In this section we describe reducing the optimization problem \rf{8no15a} conceptually.
In Section \ref{19de15b} we consider the reduction computationally.
For the remainder of this paper we fix $n$ and $q$.

We consider the natural action of $H=S_q\rtimes S_n$ on $\oR^{\CC_2}$.
If $U_1,\ldots,U_k$ form a representative set of matrices for this action, then with \rf{11de15a}
we obtain a reduction of the size of the optimization problem to polynomial size.
To make this reduction explicit in order to apply semidefinite programming, we need to express each $m_i\times m_i$ matrix $U_i\T M(y)U_i$ as an explicit
matrix in which each entry is a linear combination of the variables $y(\omega)$ for $\omega\in\Omega$ (the set of $H$-orbits of $\CC_4$).

For $\omega\in\Omega$, let $N_{\omega}$ be the $\CC_2\times\CC_2$ matrix with $0,1$ entries
satisfying
\dyz{
$(N_{\omega})_{\{\alpha,\beta\},\{\gamma,\delta\}}=1$ if and only if $\{\alpha,\beta,\gamma,\delta\}\in\omega$
}
for $\alpha,\beta,\gamma,\delta\in[q]^n$.
Then
\dyyz{
U_i\T M(y)U_i=\sum_{\omega}y(\omega)U_i\T N_{\omega}U_i.
}
So to get the reduction,
we need to obtain the matrices $U_i\T N_{\omega}U_i$ explicitly,
for each $\omega\in\Omega$ and for each $i=1,\ldots,k$.
We do this in a number of steps.

We first describe in Section \ref{16de15i}.\ref{16de15gx} a representative set for the natural
action of $S_q$ on $\oR^{q\times q}$.
From this we derive, in Section \ref{16de15i}.\ref{28de15a}, with the help of the representative set for the
action of $S_n$ on $V^{\otimes n}$ described in Section \ref{16de15j}.\ref{16de15h}, a representative set for the
action of the wreath product $H=S_q^n\rtimes S_n$ on the set
$([q]^n)^2$ of {\em ordered} pairs of words in $[q]^n$, in other words,
on $\oR^{([q]^n)^2}\cong(\oR^{q\times q})^{\otimes n}$.
From this we derive in Section \ref{16de15i}.\ref{28de15b} a
representative set for the action of $H$ on the set $\CC_2\setminus\{\emptyset\}$
of {\em unordered} pairs $\{v,w\}$ (including singleton) of words $v,w$
in $[q]^n$.
Then in Section \ref{16de15i}.\ref{28de15c} we derive a representative
set for the action of $H$ on the set $\CC_2^d\setminus\{\emptyset\}$,
where $\CC_2^d$ is the set of codes in $\CC_2$ of minimum distance at least $d$.
(So each singleton word belongs to $\CC_2^d$.)
Finally, in Section \ref{16de15i}.\ref{28de15c} we include the
empty set $\emptyset$, by an easy representation-theoretic argument.

\subsect{16de15gx}{A representative set for the action of $S_q$ on $\oR^{q\times q}$}

We now consider the natural action of $S_q$ on $\oR^{q\times q}$.
Let $e_j$ be the $j$-th unit basis vector in $\oR^q$, $J_q$ be the all-one
$q\times q$ matrix, $\bm{1}$ be the all-one column vector in $\oR^q$,
$N:=(e_1-e_2)\bm{1}\T$, and $E_{i,j}:=e_ie_j\T$.
We furthermore define the following matrices, where we consider matrics
in $\oR^{q\times q}$ as {\em columns} of the matrices $B_i$:
\dy{9no15a}{
$B_1:=[I_q,J_q-I_q]$,\\
$B_2:=[E_{1,1}-E_{2,2},
N-N\T, N+N\T-2(E_{1,1}-E_{2,2})] $,\\
$B_3:=
[E_{1,2}+E_{2,3}+E_{3,1}-E_{2,1}-E_{3,2}-E_{1,3}]$,\\
$B_4:=
[E_{1,3}-E_{3,2}+E_{2,4}-E_{4,1}+E_{3,1}-E_{2,3}+E_{4,2}-E_{1,4}]$.
}
The matrices in $\oR^{q\times q}$ will in fact be taken as elements of the dual space
$(\oR^{q\times q})^*$ (by taking the inner product), so that they are elements of the algebra
$\OO(\oR^{q\times q})$ of polynomials on the linear space $\oR^{q\times q}$.

One may check that $\{B_1,\ldots,B_4\}$ is representative for the natural action of
$S_q$ on $\oR^{q\times q}$, if $q\geq 4$.
If $q\leq 3$, we delete $B_4$, and if $q=2$ we moreover delete $B_3$ and the last column of $B_2$
(as this column is 0 if $q=2$).
Therefore, if $q\geq 4$, set $k=4$, $m_1=2$, $m_2:=3$, $m_3:=1$, and $m_4:=1$.
If $q=3$, set $k=3$, $m_1=2$, $m_2:=3$, and $m_3:=1$.
If $q=2$, set $k=2$, $m_1=2$, and $m_2:=2$.

For the remainder of this paper we fix $k$, $m_1,\ldots,m_k$, and $B_1,\ldots,B_k$.

\subsect{28de15a}{A representative set for the action of $H$ on $(\oR^{q\times q})^{\otimes n}$}

We next consider the action of $H$ on
the set $([q]^n)^2$ of {\em ordered} pairs of code words.
For that, we derive a representative set for the natural action of $H$
on $(\oR^{q\times q})^{\otimes n}\cong\oR^{([q]^n)^2}$ from the results described in
Sections \ref{16de15j}.\ref{16de15h} and \ref{16de15i}.\ref{16de15gx}.

Let {\bm{N}} be the collection of all $k$-tuples $(n_1,\ldots,n_k)$ of nonnegative integers adding up
to $n$.
For $\bm{n}=(n_1,\ldots,n_k)\in{\bm N}$, let
$\bm{\lambda\vdash n}$ mean that $\bm{\lambda}=(\lambda_1,\ldots,\lambda_k)$ with
$\lambda_i\vdash n_i$ for $i=1,\ldots,k$.
(So each $\lambda_i$ is equal to $(\lambda_{i,1},\ldots,\lambda_{i,t})$ for some $t$.)

For $\bm{\lambda}\vdash\bm{n}$ define
\dez{
W_{\bm{\lambda}}:=T_{\lambda_1,m_1}\times\cdots\times T_{\lambda_k,m_k},
}
and for $\bm{\tau}=(\tau_1,\ldots,\tau_k)\in W_{\bm{\lambda}}$ define
\de{16de15d}{
v_{\bm{\tau}}:=\bigotimes_{i=1}^ku_{\tau_i,B_i}.
}

\prop{19se15a}{
The matrix set
\dyy{19se15b}{
\{~~
[v_{\bm{\tau}}
\mid
\bm{\tau}\in W_{\bm{\lambda}}]
~~\mid
\bm{n}\in\bm{N},\bm{\lambda\vdash n}
\}
}
is representative for the action of $S_q^n\rtimes S_n$ on $(\oR^{q\times q})^{\otimes n}$.
}

\pf
Let $L_i$ denote the linear space spanned by $B_i(1),\ldots,B_i(m_i)$.
Then
\dyyz{
(\oR^{q\times q})^{\otimes n}
\stackrel{\text{by \rf{28de15ster}}}{=}
\big(\bigoplus_{i=1}^k\bigoplus_{j=1}^{m_i}\oR S_q\cdot B_i(j)\big)^{\otimes n}
=
\oR S_n\cdot\bigoplus_{\bm{n}\in\bm{N}}\bigotimes_{i=1}^k
\big(\bigoplus_{j=1}^{m_i}\oR S_q\cdot B_i(j)\big)^{\otimes n_i}
=$\\
$\dps
\oR S_n\cdot\oR S_q^{\otimes n}\cdot
\bigoplus_{\bm{n}\in\bm{N}}\bigotimes_{i=1}^kL_i^{\otimes n_i}
\stackrel{\text{by \rf{28de15ster}}}{=}
\oR H\cdot
\bigoplus_{\bm{n}\in\bm{N}}\bigotimes_{i=1}^k\bigoplus_{\lambda_i\vdash n_i}\bigoplus_{\tau_i\in T_{\lambda_i,m_i}}
\oR S_{n_i}\cdot u_{\tau_i,B_i}
=
\bigoplus_{\bm{n}\in\bm{N}}
\bigoplus_{\bm{\lambda\vdash n}}
\bigoplus_{\bm{\tau}\in W_{\bm{\lambda}}}
\oR H\cdot v_{\bm{\tau}}.
}
Now for each {\bm{n,\lambda}} and $\bm{\tau},\bm{\sigma}\in W_{\bm{\lambda}}$,
there is an $H$-isomorphism
$\oR H\cdot v_{\bm{\tau}}\to\oR H\cdot v_{\bm{\sigma}}$ bringing
$v_{\bm{\tau}}$ to $v_{\bm{\sigma}}$,
since for each $i=1,\ldots,k$, setting $H_i:=S_q^{n_i}\rtimes S_{n_i}$,
there is an $H_i$-isomorphism $\oR H_i\cdot u_{\tau_i,B_i}\to\oR H_i\cdot u_{\sigma_i,B_i}$.
Hence
\dyy{16se15c}{
\dim((\oR^{q\times q})^{\otimes n}\otimes (\oR^{q\times q})^{\otimes n})^H
\geq
\sum_{\bm{n}\in\bm{N}}
\sum_{\bm{\lambda\vdash n}}|W_{\bm{\lambda}}|^2
=
\sum_{\bm{n}\in\bm{N}}
\sum_{\bm{\lambda\vdash n}}\prod_{i=1}^k|T_{\lambda_i,m_i}|^2
=
\sum_{\bm{n}\in\bm{N}}\prod_{i=1}^k\sum_{\lambda_i\vdash n_i}|T_{\lambda_i,m_i}|^2
\stackrel{\text{by \rf{16se15a}}}{=}
\sum_{\bm{n}\in\bm{N}}\prod_{i=1}^k\dim\Sym_{n_i}(\oR^{m_i}\otimes\oR^{m_i})
=
\sum_{\bm{n}\in\bm{N}}\prod_{i=1}^k\kies{m_i^2+n_i-1}{n_i-1}
=
\kies{\sum_{i=1}^km_i^2+n-1}{n-1}
=
\dim\Sym_n(((\oR^{q\times q})\otimes (\oR^{q\times q}))^{S_q})
=
\dim((\oR^{q\times q})^{\otimes n}\otimes (\oR^{q\times q})^{\otimes n})^H
}
as $\sum_{i=1}^km_i^2=\dim(\oR^{q\times q}\otimes\oR^{q\times q})^{S_q}$.
So we have equality throughout in \rf{16se15c}, and hence each $\oR H\cdot v_{\bm{\tau}}$ is irreducible, and if
$\bm{\lambda}\neq\bm{\lambda'}$,
then for each
${\bm{\tau}}\in W_{\bm{\lambda}}$
and
${\bm{\tau'}}\in W_{\bm{\lambda'}}$,
$\oR H\cdot v_{\bm{\tau}}$ 
and
$\oR H\cdot v_{\bm{\tau'}}$ 
are not $H$-isomorphic.
\bx

\subsect{28de15b}{Unordered pairs}

We now go over from the set $([q]^n)^2$ of ordered pairs of code words to the set $\CC_2\setminus\{\emptyset\}$ of
unordered pairs (including singletons) of code words.
For this we consider the action of the group $S_2$ on
$\oR^{[q]^n\times[q]^n}\cong\oR^{([q]^n)^2}\cong(\oR^{q\times q})^{\otimes n}$,
where the nonidentity element $\sigma$ in $S_2$ acts as taking the transpose.
The actions of $S_2$ and $H$ commute.

Let $F$ be the $(\CC_2\setminus\{\emptyset\})\times ([q]^n)^2$ matrix with $0,1$ entries
satisfying
\dyz{
$F_{\{\alpha,\beta\},(\gamma,\delta)}=1$ if and only if $\{\gamma,\delta\}=\{\alpha,\beta\}$,
}
for $\alpha,\beta,\gamma,\delta\in[q]^n$.
Then the function $x\mapsto Fx$ is an $H$-isomorphism
$(\oR^{([q]^n)^2})^{S_2}\to\oR^{\CC_2\setminus\{\emptyset\}}$.

Now note that each $B_i(j)$, as matrix in $\oR^{q\times q}$, is $S_2$-invariant (i.e.,\ symmetric)
except for $B_2(2)$ and $B_3(1)$,
while $\sigma\cdot B_2(2)=-B_2(2)$ and $\sigma\cdot B_3(1)=-B_3(1)$
(as $B_2(2)$ and $B_3(1)$ are skew-symmetric).
So for any
${\bm{n}}\in\bm{N}$, $\bm{\lambda\vdash n}$, and $\bm{\tau}\in W_{\bm{\lambda}}$,
we have
\dez{
\sigma\cdot v_{\bm{\tau}}
=
(-1)^{|\tau_2^{-1}(2)|+|\tau_3^{-1}(1)|}
v_{\bm{\tau}}.
}
Therefore, let $W'_{\bm{\lambda}}$ be the set of those $\bm{\tau}\in W_{\bm{\lambda}}$
with
$|\tau_2^{-1}(2)|+|\tau_3^{-1}(1)|$ even.
Then the matrix set
\dyyz{
\{~~
[ v_{\bm{\tau}} \mid \bm{\tau}\in W'_{\bm{\lambda}} ]
~~\mid~~
\bm{n}\in\bm{N},\bm{\lambda\vdash n}\}
}
is representative for the action of $H$ on $(\oR^{([q]^n)^2})^{S_2}$.
Hence the matrix set
\dyy{21de15a}{
\{~~
[Fv_{\bm{\tau}} \mid \bm{\tau}\in W'_{\bm{\lambda}} ]
~~{\mid}~~
\bm{n}\in\bm{N}, \bm{\lambda\vdash n}\}
}
is representative for the action of $H$ on $\oR^{\CC_2\setminus\{\emptyset\}}$.

\subsect{28de15c}{Restriction to pairs of words at distance at least $d$}

Let $d\in\oZ_+$, and let
$\CC^d_2$ be the collection of elements of $\CC_2$ of minimum distance at least $d$.
Note that each singleton code word belongs to $\CC_2^d$, and that $H$ acts on $\CC^d_2$.
From \rf{21de15a} we derive a representative set for the action of $H$ on $\oR^{\CC_2^d\setminus\{\emptyset\}}$.

To see this, let for each $t\in\oZ_+$, $L_t$ be the subspace of $\oR^{\CC_2}$ spanned by the elements $e_{\{\alpha,\beta\}}$ with
$\alpha,\beta\in[q]^n$ and $d_H(\alpha,\beta)=t$.
(For any $Z\in\CC_2$, $e_Z$ denotes the unit base vector in $\oR^{\CC_2^d}$ for coordinate $Z$.)

Then for any
${\bm{n}\in\bm{N}}$, ${\bm{\lambda\vdash n}}$, and $\bm{\tau}\in W'_{\bm{\lambda}}$,
the irreducible representation $H\cdot Fv_{\bm{\tau}}$ is contained in $L_t$, where
\dez{
t:=n-|\tau_1^{-1}(1)|-|\tau_2^{-1}(1)|,
}
since $B_1(1)=I_q$ and $B_2(1)=E_{1,1}-E_{2,2}$ are the only two entries $B_i(j)$ in the $B_i$ that have nonzeros on the diagonal
of the matrix $B_i(j)$.
Let $W''_{\bm{\lambda}}$ be the set of those $\bm{\tau}$ in $W'_{\bm{\lambda}}$
with
\dez{
n-|\tau_1^{-1}(1)|-|\tau_2^{-1}(1)|\in \{0,d,d+1,\ldots,n\}.
}
Then a representative set for the action of $H$ on $\CC^d_2\setminus\{\emptyset\}$ is 
\dez{
\big\{~~
[Fv_{\bm{\tau}}
\mid
\bm{\tau}\in W''_{\bm{\lambda}}
]
~~\mid~~
\bm{n}\in\bm{N}, \bm{\lambda\vdash n}
\big\}.
}

\subsect{28de15d}{Adding $\emptyset$}

To obtain a representative set for the action of $H$ on $\CC^d_2$, note that $H$ acts trivially
on $\emptyset$.
So $e_{\emptyset}$ belongs to the $H$-isotopical component of $\oR^{\CC_2}$ that consists
of $H$-invariant elements.
Now the $H$-isotypical component of $\oR^{\CC_2\setminus\{\emptyset\}}$ that consists of
the $H$-invariant elements corresponds to the matrix in the representative set indexed by
indexed by $\bm{n}=(n,0,0,0)$ and $\bm{\lambda}=((n),(),(),())$, where $()\vdash 0$.
So to obtain a representative set for $\oR^{\CC_2}$, we just add $e_{\emptyset}$ as column to this matrix.

\sect{19de15b}{How to compute $(Fv_{\bm{\tau}})\T N_{\omega}Fv_{\bm{\sigma}}$}

We now have a reduction of the original problem to blocks with coefficients
$(Fv_{\bm{\tau}})\T N_{\omega}Fv_{\bm{\sigma}}$,
for $\bm{n}\in\bm{N}$, $\bm{\lambda\vdash n}$, $\bm{\tau},\bm{\sigma}\in W_{\bm{\lambda}}$,
and $\omega\in\Omega$.
The number and orders of these blocks are bounded by a polynomial in $n$, but computing 
these coefficients still must be reduced in time, since the order of
$F$, $v_{\bm{\tau}}$, $v_{\bm{\sigma}}$, and $N_{\omega}$ is exponential in $n$.

Fix $\bm{n}\in\bm{N}$, $\bm{\lambda\vdash n}$, and $\bm{\tau},\bm{\sigma}\in W_{\bm{\lambda}}$.
For any $\omega\in\Omega$, let $L_{\omega}:=F\T N_{\omega}F$.
So $L_{\omega}$ is a $([q]^n\times [q]^n)\times([q]^n\times [q]^n)$
matrix with 0,1 entries satisfying
\dyz{
$(L_{\omega})_{(\alpha,\beta),(\gamma,\delta)}=1$ if and only if $\{\alpha,\beta,\gamma,\delta\}\in\omega$,
}
for all $\alpha,\beta,\gamma,\delta\in[q]^n$.
By definition of $L_{\omega}$,
\dez{
(Fv_{\bm{\tau}})\T N_{\omega}Fv_{\bm{\sigma}}
=
v_{\bm{\tau}}\T L_{\omega}v_{\bm{\sigma}}.
}
So it suffices to evaluate the latter value.

Let $\Pi$ be the collection of partitions of $\{1,2,3,4\}$ into at most $q$ parts.
There is the following bijection between $\Pi$ and the
set of orbits of the action of $S_q$ on $[q]^4$.

For each word $w\in[q]^4$, let $\part(w)$ be the partition $P\in\Pi$ such that $i$ and $j$ belong to
the same class of $P$ if and only if $w_i=w_j$ (for $i,j=1,\ldots,4$).
Then two elements $v,w\in[q]^4$ belong to the same $S_q$-orbit if and only if $\part(v)=\part(w)$.
Note that $|\Pi|=8$ if $q=2$, $|\Pi|=14$ if $q=3$, and $|\Pi|=15$ if $q\geq 4$.
(In all cases, $|\Pi|=\dim(\oR^{q\times q})^{S_q}=\sum_{i=1}^km_i^2$.)

For $P\in\Pi$, let
\dez{
d_P:=\sum_{i_1,\ldots,i_4\in[q]\atop\part{i_1\cdots i_4}=P}e_{i_1}e_{i_2}\T\otimes e_{i_3}e_{i_4}\T,
}
where each $e_i$ is a unit basis column vector in $\oR^q$, so that $e_ie_j\T$ is a matrix
in $\oR^{q\times q}$.
Then $D:=\{d_P\mid P\in\Pi\}$ is a basis of $(\oR^{q\times q}\otimes\oR^{q\times q})^{S_q}$.
Let $D^*$ be the dual basis.

For any $(\alpha,\beta,\gamma,\delta)\in([q]^n)^4$, let
\dez{
\psi(\alpha,\beta,\gamma,\delta):=\prod_{i=1}^nd^*_{\part(\alpha_i\beta_i\gamma_i\delta_i)},
}
which is a degree $n$ polynomial on $(\oR^{q\times q}\otimes\oR^{q\times q})^{S_q}$.
Then $\psi(\alpha,\beta,\gamma,\delta)=\psi(\alpha',\beta',\gamma',\delta')$
if and only if $(\alpha,\beta,\gamma,\delta)$ and $(\alpha',\beta',\gamma',\delta')$
belong to the same $H$-orbit on $([q]^n)^4$.
So this gives a bijection between
the set $Q$ of degree $n$ monomials expressed in the dual basis $D^*$
and the set of $H$-orbits on $([q]^n)^4\cong ([q]^4)^n$.
The function $([q]^n)^4\to\CC_4$ with $(\alpha_1,\ldots,\alpha_4)\mapsto\{\alpha_1,\ldots,\alpha_4\}$
then gives a surjective function $\omega:Q\to\Omega\setminus\{\{\emptyset\}\}$.

For any $\mu\in Q$, define
\dez{
K_{\mu}:=\sum_{d_1,\ldots,d_n\in D\atop d_1^*\cdots d_n^*=\mu}\bigotimes_{j=1}^nd_j.
}

\lemmaz{
For each $\omega\in\Omega$:
$\dps L_{\omega}=\sum_{\mu\in Q\atop\omega(\mu)=\omega}K_{\mu}$.
}

\pf
Choose $\alpha,\beta,\gamma,\delta\in [q]^n$.
Then
\dyyz{
\sum_{\mu\in Q\atop\omega(\mu)=\omega}
(K_{\mu})_{(\alpha,\beta),(\gamma,\delta)}
=
\sum_{\mu\in Q\atop\omega(\mu)=\omega}
\sum_{P_1,\ldots,P_n\in\Pi\atop d^*_{P_1}\cdots d^*_{P_n}=\mu}
\big(\bigotimes_{i=1}^nd_{P_i}\big)_{\alpha,\beta,\gamma,\delta}
=
\sum_{\mu\in Q\atop\omega(\mu)=\omega}
\sum_{P_1,\ldots,P_n\in\Pi\atop d^*_{P_1}\cdots d^*_{P_n}=\mu}
\prod_{i=1}^n(d_{P_i})_{\alpha_i,\beta_i,\gamma_i,\delta_i}.
}
Now the latter value is 1 if
$\omega(d^*_{\part(\alpha_1\beta_1\gamma_1\delta_1)}\cdots d^*_{\part(\alpha_n\beta_n\delta_n\gamma_n)})=\omega$,
and is 0 otherwise.
So it is equal to $(L_{\omega})_{(\alpha,\beta),(\gamma,\delta)}$.
\bx

By this lemma, it suffices to compute $v_{\bm{\tau}}\T K_{\mu}v_{\bm{\sigma}}$ for each
$\mu\in Q$.
To this end, define the following degree $n$ polynomial
on $W:=(\oR^{q\times q}\otimes\oR^{q\times q})^{S_q}$:
\de{14de15d}{
p_{\bm{\tau,\sigma}}:=
\prod_{i=1}^k
\sum_{\tau_i'\sim\tau_i\atop\sigma_i'\sim\sigma_i}\sum_{c_i,c_i'\in C_{\lambda_i}}\sgn(c_ic_i')
\prod_{y\in Y(\lambda_i)}
B_i(\tau_i'c_i(y))\otimes B_i(\sigma_i'c_i'(y)).
}
This polynomial can be computed (i.e., expressed as linear combination of monomials
in $B_i(j)\otimes B_i(h)$) in time bounded by a
polynomial in $n$ (Gijswijt [5], see Appendix 1 in Section \ref{19de15b}.\ref{21de15b} below).

\lemmaz{
$\dps\sum_{\mu\in Q}(v_{\bm{\tau}}\T K_{\mu}v_{\bm{\sigma}})\mu
=
p_{\bm{\tau,\sigma}}$.
}

\pf
We can write for each $\mu\in Q$:
\dez{
v_{\bm{\tau}}\T K_{\mu}v_{\bm{\sigma}}
=
(v_{\bm{\tau}}\otimes v_{\bm{\sigma}})(K_{\mu}),
}
using the fact that
$v_{\bm{\tau}},v_{\bm{\sigma}}\in((\oR^{q\times q})^{\otimes n})^*$ and
$K_{\mu}\in(\oR^{q\times q})^{\otimes n}\otimes(\oR^{q\times q})^{\otimes n}$.
So it suffices to show
\dyy{29de15b}{
\sum_{\mu\in Q}(v_{\bm{\tau}}\otimes v_{\bm{\sigma}})(K_{\mu})\mu=p_{\bm{\tau,\sigma}}.
}

Consider any $f=f_1\cdots f_n$ with $f_j\in W^*$ for $j=1,\ldots,n$.
Then
\dyy{16de15g}{
f=\sum_{\mu\in Q}(\bigotimes_{j=1}^nf_j)(K_{\mu})\mu.
}
Indeed,
\dyyz{
\sum_{\mu\in Q}(\bigotimes_{j=1}^nf_j)(K_{\mu})\mu
=
\sum_{d_1,\ldots,d_n\in D\atop d_1^*\cdots d_n^*=\mu}(\bigotimes_{j=1}^nf_j)(\bigotimes_{j=1}^n d_j)\mu
=
\sum_{d_1,\ldots,d_n\in D}\prod_{j=1}^nf_j(d_j)d_j^*
=
\prod_{j=1}^n\sum_{d\in D}f_j(d)d^*
=
\prod_{j=1}^nf_j
=
f.
}
Applying \rf{16de15g} to each term $f$ of $p_{\bm{\tau,\sigma}}$ as given by \rf{14de15d} we obtain
\rf{29de15b}, in view of \rf{14de15f} and \rf{16de15d}.
\bx

So $v_{\bm{\tau}}\T K_{\mu}v_{\bm{\sigma}}$ can be computed by
expressing the polynomial $p_{\bm{\tau,\sigma}}$ as linear combination of monomials $\mu\in Q$,
which are products of linear functions in $D^*$.
So it suffices to express each $B_i(j)\otimes B_i(h)$
as linear function into the basis $D^*$, that is, to
calculate the numbers $(B_i(j)\otimes B_i(h))(d_P)$ for all $i=1,\ldots,k$,
$j,h=1,\ldots,m_i$, and $P\in\Pi$
 --- see Appendix 2 (Section \ref{19de15b}.\ref{19de15a} below).

We finally consider the entries in the row and column for $\emptyset$ in
the matrix associated with $\bm{\lambda}=((n),(),(),())$
(cf.\ Section \ref{16de15i}.\ref{28de15d}).
Trivially, $e_{\emptyset}\T M(x)e_{\emptyset}=(M(x))_{\emptyset,\emptyset} =x(\emptyset)$, which is set to 1 in the optimization problem.
Any $\bm{\tau}\in W_{\bm{\lambda}}$ is determined by the number $t$ of 2's in the row of the Young shape $Y((n))$.
Then
\dyyz{
v_{\bm{\tau}}=\sum_{u,w\in[q]^n\atop d_H(u,w)=t}e_{(u,w)}
\text{~~~and hence~~~}
Fv_{\bm{\tau}}=\sum_{u,w\in[q]^n\atop d_H(u,w)=t}e_{\{u,w\}}.
}
Hence, as $\emptyset\cup\{u,w\}=\{u,w\}$,
\dyyz{
e_{\emptyset}\T M(x)Fv_{\bm{\tau}}
=
\sum_{u,w\in[q]^n\atop d_H(u,w)=t}x({\{u,w\}})
=
\kies{n}{t}q^n(q-1)^ty(\omega),
}
where $\omega$ is the $H$-orbit of $\CC_4$ consisting of all pairs $\{\alpha,\beta\}$ with $d_H(\alpha,\beta)=t$.

\subsect{21de15b}{Appendix 1: Computation of $p_{\tau,\sigma}$}

For any $n,m\in\oZ_+$, $\lambda\vdash n$, and $\tau,\sigma\in T_{\lambda,m}$, define
the polynomial $p_{\tau,\sigma}\in\oR[x_{j,h}\mid j,h=1,\ldots,m]$ by
\dyyz{
p_{\tau,\sigma}(X):=
\sum_{\tau'\sim\tau\atop\sigma'\sim\sigma}\sum_{c,c'\in C_{\lambda}}\sgn(cc')
\prod_{y\in Y(\lambda)}
x_{\tau'c(y),\sigma'c'(y)},
}
for $X=(x_{j,h})_{j,h=1}^m\in\oR^{m\times m}$.

\prop{19se15d}{
Expressing $p_{\tau,\sigma}$ as a linear combination of monomials can be done in polynomial time, for fixed $m$.
}

\pf
First observe that
\dyyz{
p_{\tau,\sigma}(X)
=
|C_{\lambda}|\sum_{\tau'\sim\tau\atop\sigma'\sim\sigma}\sum_{c\in C_{\lambda}}\sgn(c)
\prod_{y\in Y_{\lambda}}x_{\tau'(y),\sigma'c(y)}
=
|C_{\lambda}|
\sum_{\tau'\sim\tau\atop\sigma'\sim\sigma}\prod_{j=1}^{\lambda_1}\det((x_{\tau'(i,j),\sigma'(i',j)})_{i,i'=1}^{\lambda^*_j}).
}
($\lambda^*$ is the dual partition of $\lambda$; that is, $\lambda^*_j$ is the height of column $j$.)

For fixed $m$, when $n$ grows, there will be several columns of $Y(\lambda)$ that are the same both in $\tau'$
and in $\sigma'$.
More precisely, for given $\tau',\sigma'$ let the `count function' $\kappa$ be defined as follows:
for $t\in\oZ_+$ and $v,w\in[m]^t$, $\kappa(v,w)$
is the number of columns $j$ of height $t$ such that
$\tau'(i,j)=v_i$ and $\sigma'(i,j)=w_i$ for all $i=1,\ldots,t$.
Then for each $i\leq h:=\height(\lambda)$ and each $s\in[m]$:
\dy{5se15d}{
$\dps\sum_{t=i}^{h}\sum_{v,w\in[m]^t\atop v_i=s}\kappa(v,w)=$ number of $s$ in row $i$ of $\tau$, and\\
$\dps\sum_{t=i}^{h}\sum_{v,w\in[m]^t\atop w_i=s}\kappa(v,w)=$ number of $s$ in row $i$ of $\sigma$.
}
For any given function $\kappa:\bigcup_{i=1}^{h}[m]^i\times[m]^i\to\oZ_+$ satisfying \rf{5se15d}, there are precisely
\de{29de15a}{
\prod_{t=1}^h\frac{(\lambda_t-\lambda_{t+1})!}{\prod_{v,w\in[m]^t}\kappa(v,w)!}
}
pairs $\tau'\sim\tau$ and $\sigma'\sim\sigma$ having count function $\kappa$ (setting $\lambda_{h+1}:=0$).
(Note that \rf{5se15d} implies $\lambda_t-\lambda_{t+1}=\sum_{v,w\in[m]^t}\kappa(v,w)$, for each $t$,
so that for each $t$, the factor in \rf{29de15a} is a Newton multinomial coefficient.)
Hence
\dyyz{
p_{\tau,\sigma}
=
|C_{\lambda}|
\sum_{\kappa}\prod_{t=1}^h(\lambda_t-\lambda_{t+1})!\prod_{v,w\in[m]^t}
\frac{\det((x_{v(i),w(i')})_{i,i'=1}^t)^{\kappa(v,w)}}{\kappa(v,w)!},
}
where $\kappa$ ranges over functions $\kappa:\bigcup_{t=1}^h([m]^t\times [m]^t)\to\oZ_+$ satisfying \rf{5se15d}.
\bx

\subsect{19de15a}{Appendix 2: Expressing $B_i(j)\otimes B_i(h)$ into $d^*_P$}

Recall that each $B_i(j)$ is a linear function on $\oR^{q\times q}$, and that each $d_P$ is an element of
$\oR^{q\times q}\otimes\oR^{q\times q}$, where $P$ belongs to the set $\Pi$ of partitions of
$\{1,\ldots,4\}$ with at most $q$ classes.
We express each $B_i(j)\otimes B_i(h)$ in the dual basis $B^*:=\{d^*_P\mid P\in\Pi\}$.
The coefficient of $d^*_P$ is obtained by evaluating $(B_i(j)\otimes B_i(h))(d_P)$.
This is routine, but we display the expressions.

For this, denote any subset $X$ of $\{1,\ldots,4\}$ by a string formed by the elements of $X$,
and denote a partition $P$ of $\{1,\ldots,4\}$  by a sequence of its classes
(for instance, $d^*_{13,2,4}$ denotes the dual variable $d^*_P$ associated with partition
$P=\{\{1,3\},\{2\},\{4\}\}$ of $\{1,2,3,4\}$).
Then:

{\small
\medskip
\noindent
$B_1(1)\otimes B_1(1)=
q d^*_{1234} +q(q-1) d^*_{12,34}$,\\
$B_1(1)\otimes B_1(2)=
q(q-1)(d^*_{123,4} + d^*_{124,3} +(q-2) d^*_{12,3,4})$,\\
$B_1(2)\otimes B_1(1)=
q(q-1)(d^*_{1,234} + d^*_{134,2} +(q-2) d^*_{1,2,34})$,\\
$B_1(2)\otimes B_1(2)=
q(q-1)(d^*_{13,24} + d^*_{14,23} +(q-2)(d^*_{13,2,4} + d^*_{14,2,3} + d^*_{1,23,4} + d^*_{1,24,3} +(q-3) d^*_{1,2,3,4}))$.

\medskip
\noindent
$B_2(1)\otimes B_2(1)=
2 d^*_{1234} -2 d^*_{12,34}$,\\
$B_2(1)\otimes B_2(2)=
2q(d^*_{123,4} - d^*_{124,3})$,\\
$B_2(1)\otimes B_2(3)=
2(q-2)(d^*_{124,3} + d^*_{123,4} -2 d^*_{12,3,4})$,\\
$B_2(2)\otimes B_2(1)=
2q(d^*_{134,2} - d^*_{1,234})$,\\
$B_2(2)\otimes B_2(2)=
2q(2d^*_{13,24} -2 d^*_{14,23} +(q-2)(d^*_{13,2,4} - d^*_{14,2,3} - d^*_{1,23,4} + d^*_{1,24,3}))$,\\
$B_2(2)\otimes B_2(3)=
2q(q-2)(d^*_{13,2,4} + d^*_{14,2,3} - d^*_{1,23,4} - d^*_{1,24,3})$,\\
$B_2(3)\otimes B_2(1)=
2(q-2)(d^*_{1,234} + d^*_{134,2} -2 d^*_{1,2,34})$,\\
$B_2(3)\otimes B_2(2)=
2q(q-2)(d^*_{13,2,4} - d^*_{14,2,3} + d^*_{1,23,4} - d^*_{1,24,3})$,\\
$B_2(3)\otimes B_2(3)=
2(q-2)(2 d^*_{13,24} +2 d^*_{14,23} +(q-4)(d^*_{13,2,4} + d^*_{14,2,3} + d^*_{1,23,4} + d^*_{1,24,3}) +4(q-3) d^*_{1,2,3,4})$.

\medskip
\noindent
$B_3(1)\otimes B_3(1)=
6(d^*_{13,24} - d^*_{14,23} - d^*_{13,2,4} + d^*_{14,2,3} + d^*_{1,23,4} - d^*_{1,24,3})$.

\medskip
\noindent
$B_4(1)\otimes B_4(1)=
8(d^*_{13,24} + d^*_{14,23} - d^*_{13,2,4} - d^*_{14,2,3} - d^*_{1,23,4} - d^*_{1,24,3}) +16 d^*_{1,2,3,4}$.

}

\bigskip
\noindent
{\em Acknowledgements.}
We are very grateful to Hans D. Mittelmann for his help in solving the larger semidefinite
programming problems.

\section*{References}\label{REF}
{\small
\begin{itemize}{}{
\setlength{\labelwidth}{4mm}
\setlength{\parsep}{0mm}
\setlength{\itemsep}{0mm}
\setlength{\leftmargin}{5mm}
\setlength{\labelsep}{1mm}
}
\item[\mbox{\rm[1]}] G.T. Bogdanova, A.E. Brouwer, S.N. Kapralov, P.R.J. \"Osterg\rondje ard, 
Error-correcting codes over an alphabet of four elements,
{\em Discrete \& Computational Geometry} 23 (2001) 333--342.

\item[\mbox{\rm[2]}] G.T. Bogdanova, P.R.J. \"Osterg\rondje ard, 
Bounds on codes over an alphabet of five elements,
{\em Discrete Mathematics} 240 (2001) 13--19.

\item[\mbox{\rm[3]}] A.E. Brouwer, 
Tables of code bounds, 2015,
see
\url{http://www.win.tue.nl/~aeb/}

\item[\mbox{\rm[4]}] P. Delsarte, 
{\em An Algebraic Approach to the Association Schemes of Coding Theory},
Philips Research Reports Supplements 1973 No. 10,
Philips Research Laboratories, Eindhoven, 1973.

\item[\mbox{\rm[5]}] D. Gijswijt, 
Block diagonalization for algebra's associated with block codes,
2014,
ArXiv \url{http://arxiv.org/abs/0910.4515}

\item[\mbox{\rm[6]}] D.C. Gijswijt, H.D. Mittelmann, A. Schrijver, 
Semidefinite code bounds based on quadruple distances,
{\em {IEEE} Transactions on Information Theory} 58 (2012) 2697--2705. 

\item[\mbox{\rm[7]}] D. Gijswijt, A. Schrijver, H. Tanaka, 
New upper bounds for nonbinary codes based on the Terwilliger algebra
and semidefinite programming,
{\em Journal of Combinatorial Theory, Series A} 113 (2006) 1719--1731.

\item[\mbox{\rm[8]}] B.E. Sagan, 
{\em The Symmetric Group: Representations, Combinatorial Algorithms, and Symmetric Functions},
Graduate Texts in Mathematics, Vol. 203, Springer, New York, 2001.

\item[\mbox{\rm[9]}] A. Schrijver, 
New code upper bounds from the Terwilliger algebra and semidefinite
programming,
{\em {IEEE} Transactions on Information Theory} 51 (2005) 2859--2866.

\end{itemize}
}

\end{document}